\documentclass[a4paper,12pt,reqno]{amsart}
\usepackage{amssymb}
\usepackage{amsmath}

\def\i{\,\lrcorner\,}
\def\a{\alpha}
\def\b{\beta}
\def\g{\gamma}

\def\vs{\vskip .6cm}

\def\la{\langle}
\def\ra{\rangle}
\def\.{\cdot}

\def\n{\nabla}
\def\l{\lambda}
\def\s{\sigma}

\def\beq{\begin{equation}}
\def\eeq{\end{equation}}
\def\bea{\begin{eqnarray*}}
\def\eea{\end{eqnarray*}}
\def\ba{\begin{array}}
\def\ea{\end{array}}

\def\f{\varphi}

\def\o{\omega}

\def\d{{\delta}}
\def\e{\varepsilon}
\def\L{\Lambda}

\def\bp{\begin{proof}}
\def\r{\end{proof}}
\def\res{\Big|}

\def\dt{\frac{\partial}{\partial t}}
\def\pa{\partial}

\def \RM{\mathbb{R}}

\def \ZM{\mathbb{Z}}

\def\es{\,\lrcorner\,}

\def\Ric{\mathrm{Ric}}

\def\be{\begin{equation}}
\def\ee{\end{equation}}
\def\Lie{{\mathcal L\,}}

\def\tr{\mathrm{tr}}

\def\gg{\mathfrak{g}}
\def\tt{\mathfrak{t}}

\newtheorem{ede}{Definition}[section]

\newtheorem{epr}[ede]{Proposition}

\newtheorem{ath}[ede]{Theorem}

\newtheorem{elem}[ede]{Lemma}

\newtheorem{ere}[ede]{Remark}

\newtheorem{ecor}[ede]{Corollary}
 
\title{Killing Vector Fields with Twistor Derivative}
\author{Andrei Moroianu}

\address{CMLS -- {\'E}cole Polytechnique \\ UMR 7640 du
  CNRS \\  91128 Palaiseau, France\\}
\email{am@math.polytechnique.fr}

\begin{document}

\begin{abstract}
Motivated by the possible characterization of Sasakian manifolds in terms of
twistor forms, we give the complete classification of compact
Riemannian manifolds carrying a Killing vector field whose covariant
derivative (viewed as a $2$--form) is a twistor form.
\vs


\end{abstract}

\maketitle

\section{Introduction}

The concept of {\em twistor forms} on Riemannian manifolds was
introduced and intensively studied by the Japanese geometers in the
50's. Some decades later, theoretical physicists became interested in
these objects, which can be used to define quadratic first integrals
of the geodesic equation, (cf. Penrose and Walker \cite{penrose}) or to
obtain symmetries of field equations (cf. \cite{be1},
\cite{be2}). More recently, a 
new impetus in this direction of research 
was given by the work of Uwe Semmelmann \cite{uwe} (see also
\cite{bms}, \cite{au}, \cite{au1}).    

Roughly speaking, a twistor form on a Riemannian manifold $M$ 
is a differential   
$p$--form $u$ such that one of the three components of its covariant
derivative $\n u$ with respect to the Levi--Civita connection vanishes
(the two other components can be identified respectively with the
differential $d u$ and codifferential $\d u$). If moreover the
codifferential $\d u$ vanishes, $u$ is called a {\em Killing
  form}. For $p=1$, twistor forms correspond to conformal vector
fields and Killing forms correspond to Killing vector fields
via the isomorphism between $T^*M$ and $TM$ induced by the metric.

Two basic examples of manifolds carrying twistor forms are the round
spheres and Sasakian manifolds, cf. \cite[Prop. 3.2 and
Prop. 3.4]{uwe}. A common feature of these examples
is the existence of Killing $1$--forms whose exterior
derivatives are twistor $2$--forms. 

Conversely, if $\xi$ is a Killing $1$--form {\em of constant length}
with twistor derivative, then it defines a Sasakian structure (see
Proposition \ref{ks} below). It is therefore natural to drop the assumption
on the length, and to address the
question of classifying all Riemannian manifolds with this
property. 

After some preliminaries on twistor forms in Section 2, we study the
behaviour of closed twistor $2$--forms with respect to the curvature
tensor in Section 3. This is used to obtain the following dichotomy in
Section 4: if $\xi$ is a Killing $1$--form with twistor exterior
derivative, then either $\xi$ satisfies a Sasaki--type equation, or
its kernel is an integrable distribution on $M$. The two possibilities
are then studied in the last four sections, where in particular new examples of
Riemannian manifolds carrying twistor $2$--forms are exhibited. A complete
classification is obtained in the compact case, cf. Theorem \ref{clas}. 

{\sc Acknowledgments.} {It is a pleasure to thank Paul Gauduchon and
  Christophe Margerin for many enlightening discussions.} 

\vs

\section{Preliminaries}

Let $(M^n,g)$ be a Riemannian manifold.
Throughout this paper vectors and $1$--forms as well as
endomorphisms of $TM$ and two times covariant tensors are identified via the
metric. In the sequel, $\{e_i\}$ will denote a local orthonormal basis
of the tangent bundle, parallel at some point. We use Einstein's
summation convention whenever subscripts appear twice.

We refer the reader to \cite{uwe} for an extensive introduction to
twistor forms. We only recall here their definition and a few basic
properties. 

\begin{ede}A $p$--form $ u $ is a {\em twistor form}
if and only if it satisfies the equation
\begin{equation}\label{twistor}
\nabla_X  u  = 
\frac{1}{  p+1} X\es d u   - 
\frac{1}{  n-p+1}  X \wedge \d  u ,
\end{equation}
for all vector fields $X$, where $du$ denotes the exterior derivative
of $u$ and $\d u$ its codifferential. If, in addition, $u$ is
co--closed $(\d u=0)$ then $u$ is said to be a {\em Killing form}.
\end{ede}

By taking one more covariant derivative in (\ref{twistor}) and summing
over an orthonormal basis $X=e_i$ we see that every twistor $p$--form satisfies
$$\n^*\n u= \frac{1}{  p+1}\d d u+\frac{1}{  n-p+1} d \d  u.$$
Taking $p=1$ and $\d u=0$ in this formula shows that
\beq\label{ff1}\n^*\n u=\frac12\Delta u,\eeq
for every Killing $1$--form $u$. For later use we also recall here 
the usual Bochner formula holding for every $1$--form $u$:
\beq\label{ff2}\Delta u=\n^*\n u+\Ric(u).\eeq

\begin{ede} \label{sa}A Sasakian structure on $M$ is a
   Killing vector field $\xi$ of constant length, such that 
\begin{equation} \label{dsa} \n^2_{X,Y}\xi=k
\big(\langle \xi,Y\rangle X-\langle X,Y\rangle  \xi\big)
, \qquad\forall\  X,Y\in TM,
\end{equation}
for some positive constant $k$.
\end{ede}

Notice that we have extended the usual definition (which assumes
$k=1$ and $\xi$ of unit length) in order to 
obtain a class of manifolds invariant through constant rescaling. 

If we denote by $u$ the $2$--form corresponding to the skew--symmetric
endomorphism $\n\xi$, then (\ref{dsa}) is equivalent to
\beq\label{dsa1} \n_Xu=k\xi\wedge X,\qquad k>0.\eeq
In particular, if $\xi$ defines a Sasakian structure, then $d\xi$ is a
closed twistor $2$--form, a fact which was noticed by U. Semmelmann
(cf. \cite[Prop. 3.4]{uwe}). As a partial converse, we have the
following characterization of Sasakian manifolds:

\begin{epr}\label{ks}
Let $\xi$ be a Killing vector field of constant length on some Riemannian
manifold such that $d\xi$ is a twistor $2$--form. Then $\xi$ is either
parallel or defines a
Sasakian structure on $M$.
\end{epr}
\bp
We may assume that $\xi$ has unit length. Let us denote by $u$ the
covariant derivative of $\xi$
\beq\label{xi1}\n_X\xi=:u(X),\qquad \forall\  X\in TM.\eeq
It is a direct consequence of the
Kostant formula that $u$ is parallel in the direction of $\xi$ (see
Section 4 for details). Since $u$ is a closed twistor form, we have 
$$\n_Xu=\frac{1}{n-1}X\wedge \d u,\qquad \forall\  X\in TM,$$
whereas for $X=\xi$ we get that $\d u$ is collinear to $\xi$. Since
$\xi$ never vanishes, there
exists some function $f$ on $M$ such that 
\beq\label{xi2}\n_X u  =f X\wedge \xi,\qquad \forall\  X\in TM.\eeq
On the other hand, $\xi$ has unit length so
$u(\xi)=0$. Differentiating this last relation with respect to some
arbitrary vector $X$ and using (\ref{xi1}) and (\ref{xi2}) yields
\beq\label{rt}u ^2(X)=fX-f\la X,\xi\ra\xi,\qquad \forall\  X\in TM,\eeq
and in particular the square norm of $u$ (as tensor) is
$$\la u,u \ra:=\la u(e_i),u(e_i)\ra= -\la u ^2(e_i), e_i\ra=(1-n)f.$$
On the other hand, (\ref{xi1}) yields for every $X\in TM$
$$\n_X(\la u,u \ra)=2f\la  X\wedge\xi,u\ra=4fu(X,\xi)=0.$$
Thus $f$ is a constant, non--positive by (\ref{rt}). If $f=0$, 
$\xi$ is parallel, otherwise $\xi$ defines a
Sasakian structure by (\ref{xi2}).
\r

\vs

\section{Closed twistor 2--forms}

In this section $(M^n,g)$ is a (not necessarily compact)
Riemannian manifold of dimension  $n>3$. We start with the following
technical result   

\begin{epr} Let $u$ be a closed twistor $2$--form, identified with a
  skew--symmetric endomorphism of $TM$. Then, for every other
  skew--symmetric endomorphism $\o$ of $TM$ one has
\beq\label{co}(n-2)(R_\o\circ u-u\circ R_\o)=(R_u\circ \o-\o\circ R_u)
+(u\circ\Ric\circ\o-\o\circ\Ric\circ u),
\eeq
where $R_\o$ is the skew--symmetric endomorphism of $TM$ defined by
$$R_\o(X):=\frac12 R_{e_j,\o(e_j)}X.$$
\end{epr}

\bp The identification between $2$--forms and skew--symmetric
endomorphisms is given by the formula
\beq\label{st}u=\frac12 e_i\wedge u(e_i).\eeq
Depending on whether $u$ is viewed as a $2$--form or as an
endomorphism, the induced action of the curvature on it reads
\beq\label{ind}R_\omega(u)=R_{\o}e_k\wedge
u(e_k)\qquad\hbox{and}\qquad
R_\omega(u)=R_\omega\circ u-u\circ R_\omega.\eeq
Let $X$ and $Y$ be vector fields on $M$ parallel at some
point. Differentiating the twistor equation satisfied by $u$
\beq\label{f}\n_Yu=\frac{1}{1-n}Y\wedge\d u\qquad\forall\  Y\in TM\eeq
in the direction of $X$ yields
\bea
\n^2_{X,Y}u&=&\frac{1}{1-n}Y\wedge\n_X\d u=\frac{1}{n-1}Y\wedge
e_j\i\n^2_{X,e_j}u\\
&=&\frac{1}{n-1}Y\wedge e_j\i R_{X,e_j}u+\frac{1}{n-1}Y\wedge
e_j\i\n^2_{e_j,X}u \\
&=&\frac{1}{n-1}Y\wedge e_j\i R_{X,e_j}u+\frac{1}{n-1}Y\wedge
e_j\i\n_{e_j}(\frac{1}{1-n}X\wedge \d u)\\
&=&\frac{1}{n-1}Y\wedge e_j\i R_{X,e_j}u-\frac{1}{(n-1)^2}Y\wedge\n_X\d
u\\
&=&\frac{1}{n-1}Y\wedge e_j\i R_{X,e_j}u+\frac{1}{n-1}\n^2_{X,Y}u,
\eea
whence
\beq\label{iu}\n^2_{X,Y}u=\frac{1}{n-2}Y\wedge e_j\i R_{X,e_j}u.\eeq
Using the first Bianchi identity we get
$$R_u(X)= \frac12 R_{e_j,u(e_j)}X=\frac12(R_{X,u(e_j)}e_j+
R_{e_j,X}u(e_j))=R_{X,u(e_j)}e_j.$$
This, together with (\ref{ind}) and (\ref{iu}) yields
\bea (n-2)\n^2_{X,Y}u&=&Y\wedge e_j\i R_{X,e_j}e_k\wedge u(e_k)\\
&=&g(e_j,R_{X,e_j}e_k)Y\wedge u(e_k)-Y\wedge R_{X,u(e_k)}e_k\\
&=&-Y\wedge u(\Ric(X))-Y\wedge R_u(X).
\eea
After skew--symmetrizing in $X$ and $Y$ we get
$$(n-2)R_{X,Y}u=X\wedge u(\Ric(Y))+X\wedge R_u(Y)-Y\wedge
u(\Ric(X))-Y\wedge R_u(X).$$
Let now $\o$ be some skew--symmetric endomorphism of $TM$. We take $X=e_i$,
$Y=\o(e_i)$ in the previous equation and sum over $i$ to obtain:
\bea (n-2)R_\o(u)&=&\frac12\big(e_i\wedge u(\Ric(\o(e_i)))+e_i\wedge
R_u(\o(e_i))\\
&&-\o(e_i)\wedge 
u(\Ric(e_i))-\o(e_i)\wedge R_u(e_i)\big)\\
&=&e_i\wedge u(\Ric(\o(e_i)))+e_i\wedge
R_u(\o(e_i))\\
&=&(u\circ\Ric\circ\o-\o\circ\Ric\circ u)+(R_u\circ \o-\o\circ R_u),
\eea
taking into account that for every endomorphism $A$ of $TM$, the
$2$--form $e_i\wedge A(e_i)$ corresponds to the skew--symmetric
endomorphism $A-^t\!\!A$ of $TM$.
\r

\begin{ecor} \label{cor} 
If $u$ is a closed twistor $2$--form, the square of the
  endomorphism corresponding to $u$ commutes with the Ricci tensor:
$$u^2\circ \Ric=\Ric\circ u ^2.$$
\end{ecor}

\bp Taking $\o=u$ in (\ref{co}) yields
\beq\label{ru}(n-3)(R_u\circ u-u\circ R_u)=0,\eeq
so $u$ and $R_u$ commute (as we assumed $n>3$). We then have
\bea 0&=& (n-2)\tr\big(u\circ(R_\o\circ u-u\circ R_\o)\big)\\
&\stackrel{(\ref{co})}{=}&\tr\big(u\circ R_u\circ\o-u\circ\o\circ R_u+
u^2\circ\Ric\circ\o-u\circ\o\circ\Ric\circ u\big)\\
&\stackrel{(\ref{ru})}{=}&\tr(u ^2\circ\Ric\circ\o-\Ric\circ u ^2\circ\o)\\
&=&-\la \o,u ^2\circ\Ric-\Ric\circ u ^2\ra .
\eea
Since $u ^2\circ\Ric-\Ric\circ u ^2$ is skew--symmetric and the
equality above holds for every skew--symmetric endomorphism $\o$, the
corollary follows.
\r

\vs

\section{Killing vector fields with twistor derivative}

We will use the general results above in the particular setting which
interests us. No compactness assumption will be needed in this
section.

Let $\xi$ be a Killing vector field on $M$, and denote by $u$ its
covariant derivative: 
\beq\label{xi}\n_X\xi=:u(X).\eeq
By definition, $u$ is a skew--symmetric tensor,
which can be identified with $\frac12 d\xi$. Taking the covariant
derivative in (\ref{xi}) yields
\beq\label{nxi}\n^2_{X,Y}\xi=(\n_X u)(Y).\eeq
This equation, together with the Kostant formula 
\beq\label{nxi2}\n^2_{X,Y}\xi=R_{X,\xi}Y\eeq
(which holds for every Killing vector field $\xi$) shows that 
\beq\label{nf}\n_\xi u=0.\eeq
Suppose now, and throughout the remaining part of this article, that
the covariant derivative 
$u$ of $\xi$ is a twistor $2$--form. Notice that in contrast to 
Proposition \ref{ks}, we no longer assume the length of $\xi$ to be
constant. Taking $Y=\xi$ in (\ref{f}) and
using (\ref{nf}) yields 
\beq\label{col}\xi\wedge\d u=0,\eeq
so $\d u$ and $\xi$ are collinear. 
We denote by $f$ the function defined on the support of
$\xi$ satisfying $(1-n)\d u=f\xi$ (this normalization
turns out to be the most convenient one in the computations below). On
the support of $\xi$ the twistor equation (\ref{f}) then reads
\beq\label{ff}\n_Xu=fX\wedge\xi,\qquad\forall\  X\in TM.\eeq
Recall now the formula  
$$(n-2)\n^2_{X,Y}u=-Y\wedge u(\Ric(X))-Y\wedge R_u(X)$$
obtained in the previous section. 
We take the inner product with $Y$ in this formula
and sum over an orthonormal basis $Y=e_i$ to obtain:
$$-(n-2)\n_X\d u=-(n-1)(u(\Ric(X))+R_u(X)).$$
Taking the scalar product with some vector $Y$ 
in this equation and symmetrizing
the result yields
$$-\frac{n-2}{n-1}(\la\n_X\d u,Y\ra+\la\n_Y\d u,X\ra)=
\la\Ric(u(X)),Y\ra+\la\Ric(u(Y)),X\ra.$$
If we replace $Y$ by $u(Y)$ in this last equation and use Corollary
\ref{cor}, we see that the expression 
$$\la\n_X\d u,u(Y)\ra+\la\n_{u(Y)}\d u,X\ra$$
is symmetric in $X$ and $Y$, {\em i.e.}
\beq\label{sy}\la\n_X\d u,u(Y)\ra+\la\n_{u(Y)}\d u,X\ra=
\la\n_Y\d u,u(X)\ra+\la\n_{u(X)}\d u,Y\ra.
\eeq
A straightforward calculation taking (\ref{ff}) and (\ref{sy}) into
account yields
\beq\label{sym}u(\xi)\wedge df+u(df)\wedge\xi=0.\eeq
On the other hand we have $u(\xi)=\n_\xi\xi=-\frac12d(|\xi|^2)$ and
$$ X(|u|^2)=2\la\n_Xu,u\ra=2f\la X\wedge\xi,u\ra=-2f\la X,u(\xi)\ra,
$$
whence
\beq\label{gy}d(|u|^2)=-2fu(\xi)=fd(|\xi|^2).\eeq
Notice that the norm $|u|$ used here is the norm of $u$ as $2$--form,
and differs by a factor $\sqrt2$ from the norm of $u$ as tensor. More
explicitly, $|u|^2=\frac12 \la u(e_i), u(e_i)\ra$.
Taking the exterior derivative in (\ref{gy}) yields 
\beq\label{dfx}0=df\wedge d(|\xi|^2)=-2df\wedge u(\xi),\eeq
which, together with (\ref{sym}) leads to
\beq\label{we}u(df)\wedge\xi=0.\eeq
The main goal of this section is to show the following

\begin{epr}\label{prop} 
Either  $f$ is constant on $M$, or $u$ has rank $2$ on $M$ and 
$\xi\wedge u=0$. 
\end{epr}

\bp
Suppose that $f$ is non--constant. Since the support of $\xi$ (say
$M_0$) is a
dense open subset of $M$, there exists a non--empty connected open
subset $U$ of $M_0$
where $df$ does not vanish. We restrict to $U$ for the
computations below. First, (\ref{dfx}) shows that $u(\xi)$ is
collinear to $df$, which, together with (\ref{we}) implies that 
\beq\label{al}u^2(\xi)=\a\xi,\eeq
for some function $\a$ defined on $U$. 

Differentiating this relation with respect to some vector $X$ and
using (\ref{xi}) and (\ref{ff}) yields
$$(X\wedge f\xi)(u(\xi))+u((X\wedge f\xi)(\xi))+u ^3(X)=\a
u(X)+X(\a)\xi,$$
or equivalently
\beq\label{jn}u^3(X)-(f|\xi|^2+\a)u(X)=(X(\a)-f\la
X,u(\xi)\ra)\xi-f\la
X,\xi\ra u(\xi).\eeq
In terms of endomorphisms of $TM$, identified with $(2,0)$--tensors,
(\ref{jn}) becomes
$$u^3-(f|\xi|^2+\a)u=(d\a-fu(\xi))\otimes\xi-f\xi\otimes u(\xi).$$
The left hand side of this relation is clearly skew--symmetric. 
The symmetric part of the right hand side thus vanishes:
$(d\a-2fu(\xi))\odot\xi=0$, whence $d\a=2fu(\xi))$ on $U$.  
Using (\ref{gy}) we get $d\a=-d(|u|^2)$, so
\beq\label{alp}\a=-|u|^2+c\eeq
for some constant $c$. We now use (\ref{jn}) in order to compute the
trace of the symmetric endomorphism $u^2$ on $T_xM$ for some $x\in U$.
It is clear that $\xi$ and
$u(\xi)$ are linearly independent eigenvectors of $u^2$ with
eigenvalue $\a$. Let $V$ denote the orthogonal complement of
$\{\xi,u(\xi)\}$ in $T_xM$. For $X\in V$, (\ref{jn}) becomes
$u^3(X)-(f|\xi|^2+\a)u(X)=0$, so the minimal polynomial of the
endomorphism $u|_V$ divides the degree $2$ polynomial
$\l(\l-(f|\xi|^2+\a))$. Thus $u^2$ has
at most $2$ different eigenvalues on $V$:  $ f|\xi|^2+\a$ and $0$,
with multiplicities denoted by $k$ and $n-k-2$ respectively. We obtain:
$$-2|u|^2=\tr(u^2)=2\a+k(f|\xi|^2+\a)
\stackrel{(\ref{alp})}{=}2c-2|u|^2+k(f|\xi|^2+\a),$$
showing that either $ f|\xi|^2+\a $ is constant or $k=0$. In the first
case we obtain by taking the exterior derivative 
\bea 0&=&d(f|\xi|^2+\a)=|\xi|^2df+fd(|\xi|^2)+d\a\\
&\stackrel{(\ref{alp})}{=}&
|\xi|^2df+fd(|\xi|^2)-d(|u|^2)\stackrel{(\ref{gy})}{=}|\xi|^2df.
\eea
This shows that $f$ is constant on $U$, contradicting the definition
of $U$. 

We therefore get $k=0$. This means that the restriction of $u$ to the
distribution $V$ vanishes, so  
\beq\label{rank2}\frac12 d\xi=u=\frac{\xi\wedge u(\xi)}{|\xi|^2}\eeq
on $U$. In particular we get 
\beq\label{u}\xi\wedge u=0\qquad\hbox{and}\qquad u\wedge u=0\ \ 
\hbox{on}\ \  U.\eeq
It remains to show that the equation $\xi\wedge u=0$ holds on
the entire manifold $M$, not only on the (possibly small) open set
$U$. This is a consequence of the following remark. The covariant
derivatives of the $3$--form $\xi\wedge u$ and of the $4$--form
$u\wedge u$ can be computed at every point of $M_0$ using (\ref{xi})
and (\ref{ff}): 
$$\n_X(\xi\wedge u)=u(X)\wedge u+\xi\wedge (fX\wedge \xi)=\frac12 X\i
(u\wedge u)$$
$$\n_X(u\wedge u)=2fX\wedge \xi\wedge u.$$
This can be interpreted by saying that the section $(\xi\wedge
u,u\wedge u)$ of $\L^3M_0\oplus\L^4M_0$ is parallel with respect to the
covariant derivative $D$ on this bundle defined by
$$D_X(\s,\tau)=(\n_X\s-\frac12 X\i \tau, \n_X\tau-2fX\wedge \s).$$
Since a parallel section which vanishes at some point is identically
zero, (\ref{u}) implies that $\xi\wedge u$ vanishes
identically on $M_0$, thus on $M$ because $M_0$ is dense in $M$.
\r

Most of the remaining part of this paper is devoted to the study of
the two possibilities given by the above proposition.

\vs

\section{The case where $f$ is constant}

In this section we consider the case where the function
$f$ defined on the support of $\xi$ is constant and we assume that $M$
is compact. We then have

\begin{ath} If the covariant derivative $u:=\n\xi$ of a non--parallel Killing
  vector field $\xi$ on $M$ satisfies 
\beq\label{nu}\n_Xu=c\xi\wedge X\eeq
for some constant $c$, then either $\xi$ defines a Sasakian structure
on $M$, or $M$ is a space form.
\end{ath}

\bp For the reader's convenience we provide here a proof of this
rather standard fact. We start by determining the
sign of the constant $c$. From (\ref{nxi}), (\ref{nxi2}) and
(\ref{nu}) we obtain
$$R_{X,\xi}Y=\n^2_{X,Y}\xi=(\n_Xu)(Y)=(c\xi \wedge X)(Y)=c(\la \xi,Y\ra X-
\la X,Y\ra \xi).$$
Taking the trace over $X$ and $Y$ in this formula yields
$$\Ric(\xi)=-R_{e_i,\xi}e_i=(n-1)c\xi.$$
Now, the two Weitzenb\"ock formulas (\ref{ff1}) and (\ref{ff2})
applied to the Killing $1$--form $\xi$ read 
$$\n^*\n\xi=\frac12\d d\xi=\frac12  \Delta\xi\qquad\hbox{and}\qquad 
\Delta\xi=\n^*\n\xi+\Ric(\xi).$$
Thus $\Ric(\xi)=\n^*\n\xi$ so taking the scalar product with $\xi$ and
integrating over $M$ yields
$$(n-1)c|\xi|^2_{L^2}=|\n\xi|^2_{L^2}.$$
This shows that $c$ is non--negative, and $c=0$ if and only if  $\xi$ is parallel,
a case which is not of interest for us. By rescaling the metric on $M$
if necessary, we can therefore assume that $c=1$, {\em i.e.} $\xi$
satisfies the Sasakian condition (\ref{dsa})
$$\n^2_{X,Y}\xi=\la\xi,Y\ra X-\la X,Y\ra\xi.$$
If the norm of $\xi$ is constant, we are in the presence of a Sasakian
structure by Definition \ref{sa}. 

Suppose that $\l:=|\xi|^2$ is non--constant. Then the function $\l$ is
a characteristic function of the round sphere. More precisely, the
second covariant derivative of the $1$--form $d\l$ can be computed as
follows. Using the relation $\n_X\xi=u(X)$ we first get $d\l=-2u(\xi)$,
therefore (\ref{nu}) gives
$$\n_Y d\l=-2(\xi\wedge Y)(\xi)-2u ^2(Y).$$
By taking another covariant derivative with respect to some
vector $X$ (at a point where $Y$ is
assumed to be parallel) we obtain after a straightforward calculation
$$\n^2_{X,Y} d\l+2X(\l)Y+Y(\l)X+d\l \la X,Y\ra=0.$$
A classical result of Tanno and Gallot (cf. \cite{tan} or 
\cite[Corollary 3.3]{gal}), shows that if $d\l$ does not vanish identically, 
the sectional curvature of $M$ has to be constant, so $M$ is a finite
quotient of the round sphere.
\r

We end up this section by remarking that conversely, every Killing vector
field on the round sphere (and all the more on its quotients) satisfies 
(\ref{ru}). This follows for instance from \cite[Prop. 3.2]{uwe}. 
The main idea is
that the space of Killing $1$--forms (respectively of closed twistor
$2$--forms) on the sphere coincides with the
eigenspace for the least eigenvalue of the Laplace operator on co--closed
$1$--forms (respectively on closed $2$--forms), and the exterior differential 
defines an isomorphism between these two spaces.

\vs

\section{The case where $\xi\wedge u=0$}

From now on we suppose that the function $f$ defined by (\ref{ff})
is non--constant. By Proposition \ref{prop} the $3$--form $\xi\wedge
d\xi$ vanishes on $M$, thus the distribution orthogonal to $\xi$
(defined on the support of $\xi$) is integrable. 
We start by a local study of the metric, at points where $\xi$ does not vanish.

\begin{epr} \label{prl} Around every point in the support of $\xi$,
  the manifold 
  $M$ is locally isometric to a warped product $I\times_\l N$ of an
  open interval $I$ and a $(n-1)$--dimensional manifold $N$ such that the
  differential of the warping function $\l$ is a twistor $1$--form on $N$.
\end{epr}

\bp
By the integrability theorem of Frobenius, $M$ can be written locally
as a product $I\times N$ 
where $\xi=\dt$ and $N$ is a local leaf tangent to the distribution $\xi
^{\perp}$. The metric $g$ can be written
$$g=\l^2 dt^2+h_t$$
for some positive function $\l$ on $I\times N$ and some family of
Riemannian metrics 
$h_t$ on $N$. Of course, the fact that $\xi=\dt$ is Killing just means
that $\l$ and $h_t$ do not depend on $t$, {\em i.e.} $g=\l^2
dt^2+h$ is a warped product. The $1$--form $\zeta$, metric dual to
$\xi$, is just $\l^2 dt$, 
so $u=\frac 12 d\zeta = \l d\l\wedge dt$. We now express the
fact that $u$ is a twistor form on $M$ in terms of the new data
$(\l,h)$. Let $X$ denote a generic vector field on $N$, identified
with the vector field on $M$ projecting over it. 
Similarly, we will identify $1$--forms on $N$ with their pull--back on $M$. 
Since the projection $M\to N$ is a Riemannian submersion, these
identifications are compatible with the metric isomorphisms between
vectors and $1$--forms. 

The O'Neill formulas (cf. \cite[p. 206]{on}) followed by a
straightforward computation give
$$\n_\dt u=0\qquad\hbox{and}\qquad \n_X u=\l\n_Xd\l\wedge dt,\qquad\forall\ 
X\in TN,$$
where we denoted by the same symbol $\n$ the covariant derivative of
the Levi--Civita connection of $h$ on $N$. 
Taking the inner product with $X$ in the second equation and summing
over an orthonormal basis of $N$ yields $\d^M u=\l \Delta
^N \l dt$, so $u$ is a twistor form if and only if 
$$\n_Xd\l=-\frac1{n-1}X\Delta ^N\l,\qquad \forall\  X\in TN$$
which just means that $d\l$ is a twistor $1$--form on $N$. 
\r

We can express the above property of $d\l$ by the fact that its
metric dual is a {\em gradient conformal vector 
field} on $N$. These objects were intensively studied in the 70's by
several authors. In particular Bourguignon \cite{bo} has shown that a
compact manifold carrying a gradient conformal vector
field is conformally equivalent to the round sphere. The
converse of this result does not hold ({\em i.e.} not every
conformally flat metric on the sphere carries gradient conformal
vector fields, cf. Remark \ref{re} below). We study this notion in
greater detail in the next section. 
\vs

\section{Gradient conformal vector fields}

\begin{ede} A {\em gradient conformal vector field} $($denoted for convenience
  GCVF in the remaining part of this paper$)$ on a connected
  Riemannian manifold $(M^n,g)$
  is a conformal vector field $X$ whose dual $1$--form is exact:
  $X=d\l$. The function $\l$ $($defined up to a constant$)$ is called the
  {\em primitive} of $X$. 
\end{ede}

Let $X$ be a GCVF. Since $X$ is a gradient vector field, its covariant
derivative is a symmetric endomorphism, and the fact that $X$ is conformal
just means that the trace--free symmetric part of $\n X$ vanishes. Thus
$X$ satisfies the equation 
\beq\label{gcvf}\n_YX=\a Y,\qquad\forall\  Y\in TM\eeq
where $\a=-\frac{\d X}{n}$. In particular, $\Lie_X g=2\a g$.

In the neighbourhood of every point where $X$ is non--zero, the metric 
$g$ can be written 
\beq\label{f1} g=\psi(t)(dt^2+h)\eeq
for some positive function
$\psi$. Conversely, if $g$ can be written in this form, then 
$\frac{\partial}{\partial t}$ is a GCVF
whose primitive is $\Psi$ (the primitive of $\psi$ in the usual
sense).

We thus see that the existence of a GCVF does not impose hard
restrictions on the metric in general. Remarkably, if the GCVF has
zeros, the situation is much more rigid:

\begin{epr} \label{prds} Let $X$ be a GCVF on a Riemannian manifold
  $(M^n,g)$ vanishing at some $x\in M$. Then there exists an open
  neighbourhood of $x$ in $M$ on which the metric can be expressed 
in polar coordinates 
\beq \label{f2}g=ds^2+\g ^2(s){g_{S^{n-1}}},\eeq
where ${g_{S^{n-1}}}$ denotes the
canonical round metric on $S^{n-1}$ and $\g$ is some
  positive function $\g:(0,\e)\to\RM^+$. The norm of $X$ in these coordinates
is a scalar multiple of $\g$: 
\beq\label{bg}|X|=c\g.\eeq
\end{epr}

Notice that the metric defined by (\ref{f2}) is in particular of type
(\ref{f1}), as shown by the change of variable
$s(t):=\int_0^t\sqrt{\psi(r)}dr$. 

\bp Let $\tau$ be the unit tangent vector field along geodesics passing
trough $x$. From (\cite[Lemma 4]{bo}) we have that $X$ is everywhere
collinear to $\tau$. Using the Gauss Lemma, we know that the metric $g$
can be expressed as $g=ds^2+h_s$ in geodesic coordinates on some
neighbourhood $U$ of $x$,
where $h_s$ is a family of metrics on $S^{n-1}$ (of course,
$\tau=\frac{\pa}{\pa s}$ in these coordinates). Since $x$ is an
isolated zero of $X$ (cf. \cite[Corollary 1]{bo}), the norm of $X$ is a
smooth function $|X|=\b$ defined on $U-\{x\}$, and $X=\b \tau$. We
then compute
\bea \dot{h_s}&=&\Lie_\tau g=\b^{-1} \Lie_Xg+2d(\b^{-1})\odot X^\flat=
2\a\b^{-1} g-2\frac{d\b}{\b}\odot ds \\
&=& 2\a\b^{-1} h_s+2\a\b^{-1} ds^2-2\frac{d\b}{\b}\odot ds
\eea
By identification of the corresponding terms in the above equality we
obtain the differential system
$$\begin{cases}{d\b}=\a ds\\
\dot{h_s}=2\a\b^{-1} h_s
\end{cases}$$
The first equation shows that $\b$ only depends on $s$:
$\b=\b(s)=\int_0^s\a(t)dt$. The second equation yields 
\beq\label{ht}h_s=\b^2(s) h\eeq 
for some metric $h$ on $S^{n-1}$. 

We claim that $h$ is (up to a scalar multiple) the canonical round
metric on the sphere. To see this, we need to understand the family of
metrics $h_s$ on $S^{n-1}$. We identify $(T_xM,g)$ with $(\RM^n,eucl)$ and
$S^{n-1}$ is viewed as the unit sphere in $T_xM$. If
$V$ is a tangent vector to $S^{n-1}$ at some $v\in T_xM$, then
$h_s(V,V)$ is the square norm with respect to $g$ of the image of $V$
by the homothety of ratio $s$ followed by the differential at $v$ of the
exponential map $\exp_x$. In other words, 
$$h_s=s^2(\exp_x)^*(g)\res_{T_vS^{n-1}}.$$
Since the differential at the origin of the exponential map is the
identity, we get 
$$\lim_{s\to 0}\frac{h_s}{s^2}={g_{S^{n-1}}}.$$
Using this together with (\ref{ht}) shows that $\lim_{s\to
  0}\frac{\b(s)}{s}$ is a positive real number denoted by $c$ and
$c^2h={g_{S^{n-1}}}$. 

We thus have proved that $g=ds^2+\gamma ^2(s){g_{S^{n-1}}}$, where
$\g=\frac{\b}{c}$. 
\r

\vs

\section{The classification}

We turn our attention back to the original question. Recall that $\xi$
is a non--parallel Killing vector field on $(M^n,g)$ such that
$\xi\wedge d\xi=0$ and $d\xi$ is a twistor form. We distinguish
two cases, depending on whether $\xi$ vanishes or not on $M$.

{\sc Case I.} The vector field $\xi$ has no zero on $M$. The
distribution orthogonal to 
$\xi$ is then globally well--defined and integrable, its maximal leaves turn
out to be compact and can be used in order
to obtain a dimensional reduction of our problem. 

\begin{ede} Let $N$ be a Riemannian manifold, let $\l$ be a positive smooth
  function on $N$ and let $\f$ be an isometry of $N$ preserving $\l$
  $($that is, $\l\circ\f=\l)$. The quotient of the warped product $\RM\times_\l
  N$ by the free $\ZM$--action generated by $(t,x)\mapsto (t+1,\f(x))$
  is called the {\em warped mapping torus} of $\f$ with respect to
  $\lambda$ and is denoted by $N_{\l,\f}$.
\end{ede}  

\begin{epr} A compact Riemannian manifold $(M^n,g)$ carries a nowhere 
  vanishing Killing vector field $\xi$ as above if and only if it is
  isometric to a warped mapping torus $N_{\l,\f}$ where $(N^{n-1},h)$
  is a compact Riemannian manifold carrying a GCVF with primitive $\l$
  and $\f$ is an isometry of $N$ preserving $\l$. 
\end{epr}

\bp The ``if" part follows directly from the local statement given by
Proposition \ref{prl}. Suppose, conversely, that $(M,\xi)$ satisfy the
conditions above. We
denote by $\f_t$ the flow of $\xi$ and by $N_x$ the maximal leaf of 
of the integrable distribution $\xi ^\perp$. Clearly $\f_t$ maps $N_x$
isometrically over $N_{\f_t(x)}$. We claim that this
action of $\RM$ on the space of leaves of $\xi ^\perp$ is
transitive. Let $x\in M$ be an arbitrary point of $M$ and denote 
$$M_x:=\bigcup_{t\in \RM}N_{\f_t(x)}.$$
For every $y\in M_x$ we define a map $\psi:(-\e,\e)\times N_y\to M$ by
$$\psi(t,z):=\f_t(z).$$
The differential of $\psi$ at $(0,y)$ is clearly invertible, thus the
inverse function theorem ensures that the image of $\psi$ contains an
open neighbourhood of $y$ in $M$. On the other hand $M_x$ contains the
image of $\psi$ by construction, therefore $M_x$ contains an open
neighbourhood of $y$. Thus $M_x$ is open. 
For any $x,y\in M$ one
either has $M_x=M_y$ or $M_x\cap M_y=\emptyset$. Thus $M$ is a
disjoint union of open sets 
$$M=\bigcup_{x\in M}M_x$$
so by connectedness we get $M_x=M$ for all $x$.

Since the norm of $\xi$ is constant along its flow, we deduce that
$|\xi|$ attains its maximum and its minimum on each integral leaf $N_x$.
By the main theorem in \cite{bo}, each leaf is conformally
diffeomorphic to the round sphere, so in particular it is
compact. Reeb's stability theorem then ensures that the space of
leaves is a compact $1$--dimensional manifold $S$ and the natural projection 
$M\to S$ is a fibration. Hence $S$ is connected, {\em i.e.} $S\cong
S^1$. On the other hand we have a group action of $\RM$ on $S$ given
by $t(N_x):=N_{\f_t(x)}$ and $S$ is the quotient of $\RM$ by the
isotropy group of some point. Since $S$ is a manifold, this isotropy
group has to be discrete, therefore is generated by some $t_0\in
\RM$. Then clearly $M$ can be identified with the warped mapping torus 
$N_{\l,\f}$, where $N:=N_x$, $\f:=\f_{t_0}$ and the warping function
$\l$ is the restriction to $N$ of $|\xi|$.
\r

\begin{ere} \label{re} A compact Riemannian manifold admitting
  gradient conformal 
  vector fields is completely classified by one single smooth function defined
  on some closed interval and satisfying some boundary
  conditions. More precisely, such a manifold is isometric to the
  Riemannian completion of a
  cylinder $(0,l)\times S^{n-2}$ with the metric
  $dt^2+f(t)g_{S^{n-2}}$, where $f:(0,l)\to \RM^+$ is smooth and
  satisfies the boundary conditions
\beq \label{bou1}
f(t)=t^2(1+t^2a(t^2))\quad\hbox{and}\quad f(l-t)=t^2(1+t^2b(t^2)),\quad
\forall\  |t|<\e,
\eeq
 for some smooth functions $a,b:(-\e,\e)\to\RM^+$. 
\end{ere}
The proof is very
 similar to that of Theorem \ref{t8} below and will thus be omitted.

{\sc Case II.} The vector field $\xi$ has zeros on $M$. The study of
this situation is 
more involved since the distribution orthogonal to $\xi$ is no longer globally
defined. On the other hand one can prove that the orbits of $\xi$
are always closed in this case, which turns out to be crucial for the
classification. This follows from a more general statement:

\begin{epr} Let $M$ be a compact Riemannian manifold
  and let $\xi$ be a Killing vector field on $M$. If the covariant
  derivative of $\xi$ has rank $2$ $($as skew--symmetric endomorphism$)$ at
  some point $x\in M$ where $\xi$ vanishes, then $\xi$ is induced by an 
  isometric $S^1$--action on $M$, and in particular its orbits are closed. 
\end{epr}

\bp Let $Z$ denote the set of points where $\xi$ vanishes, and let $Z_0$ the
connected component of $Z$ containing $x$. It is well--known that $Z_0$
is a totally geodesic submanifold of $M$ of codimension 2 (equal to
the rank of $\nabla\xi$). Moreover, at each point of $Z_0$,
$\nabla\xi$ vanishes on all vectors tangent to $Z_0$.

Since $M$ is compact, its isometry group $G$ is also compact. The
Killing vector field $\xi$ defines an element $X$ of the Lie algebra $\gg$
of $G$. The exponential map of $G$ send the line $\RM X$ onto a (not
necessarily closed) Abelian subgroup of $G$. Let $T$ be the closure of
this subgroup and denote by $\tt$ its Lie algebra. $T$ is clearly a
compact torus. We claim that $T$ is 
actually a circle. If this were not the case, one could find an
element $Y\in\tt$ defining a Killing vector field $\zeta$ on $M$ 
non--collinear to $\xi$. Let $y$ be some point in $Z_0$. Since by
definition $\xi_y=0$ we get 
$\exp(tX)\.y=y$ for all $t\in\RM$, whence $g\.y=y$ for all $g\in T$,
thus showing that $\zeta$ vanishes on $Z_0$. Since the space of
skew--symmetric endomorphisms of $T_xM$ vanishing on $T_xZ_0$ is
one--dimensional, we deduce that $(\n\zeta)_x$ is proportional to
$(\n\xi)_x$. Finally, since a Killing vector field is determined by
its 1--jet at some point, and $\xi_x=\zeta_x=0$, we deduce that
$\zeta$ is collinear to $\xi$, a contradiction.

Therefore $T$ is a circle acting isometrically on $M$ and $\xi$ is the
Killing vector field induced by this action.
\r

Let $M_0$ denote as before the set of points where $\xi$ does not
vanish. The integrable distribution $\xi ^\perp$ is well--defined
along $M_0$ and $T$ acts freely and transitively on its maximal
integral leaves. If $(N,h)$ denotes such a maximal integral leaf,
Proposition \ref{prl} shows that $M_0$ is isometric to the warped
product $S^1\times_\l N$, $g=\l^2 d\theta ^2+h$, where $\l$ is a
positive function on $N$ 
whose gradient is a conformal vector field $X$. Since $\l$ is the
restriction of the continuous function $|\xi|$ on $M$, it attains its
maximum at some $x\in N$. Of course, $X$ vanishes at $x$.

We thus may apply Proposition \ref{prds} to the gradient conformal
vector field $X$ on $N$. The metric on $N$ can be
written $h=ds^2+\g ^2(s){g_{S^{n-2}}}$ on some neighbourhood of
$x$. The length  
of $X$, which by (\ref{f2}) is equal to $c\gamma(s)$, only depends on the
distance to $x$. Assume that $X$ vanishes at some point $y:=\exp_x(tV)$
(where $V$ is a unit vector in $T_xN$). Then it vanishes on the whole
geodesic sphere of radius $t$. On the other hand $X$ has only isolated
zeros, so the geodesic sphere $S(x,t)$ is reduced to $y$. This would
imply that $N$ is compact, homeomorphic to $S^{n-1}$, so
$M_0=S^1\times N$ is compact, too. On 
the other hand $M_0$ is open, so by connectedness $M_0=M$,
contradicting the fact that $\xi$ has zeros on $M$. 

This proves that $x$ is the unique zero of $X$ on $N$. In fact we can now say
much more about the global geometry of $M$. Recall that $M$ is the
disjoint union of $M_0$ and $Z$, where $Z$, the nodal set of $\xi$, is
a codimension 2 submanifold and $M_0=N\times S^1$ is endowed with a
warped product metric. In order to state the global result we need the
following

\begin{ede} Let $l>0$ be a positive real number and let
  $\gamma, \l:(0,l)\to \RM^+$ be two smooth functions satisfying the following
  boundary conditions: 
\beq\label{equiv1}\begin{cases}\lim_{s\to 0}\gamma(s)=0,\qquad
\lim_{s\to l}\gamma(s)>0\\
\lim_{s\to 0}\l(s)>0,\qquad
\lim_{s\to l}\l(s)=0
\end{cases}\eeq

We view the sphere $S^n$ as the {\em topological join} of $S^{n-2}$ and
$S^1$, obtained from
$[0,l]\times S^{n-2}\times S^1$ by shrinking
$\{0\}\times S^{n-2}\times S^1$ to $\{\hbox{point}\}\times S^1$ 
and by shrinking
$\{l\}\times S^{n-2}\times S^1$
to $\{\hbox {point}\}\times S^{n-2}$.

Then $S^n$, endowed with the Riemannian metric
$$g=ds^2+\gamma ^2(s){g_{S^{n-2}}}+\l^2(s)d\theta ^2$$
defined on its open submanifold $(0,l)\times S^{n-2}\times S^1$ 
is called the {\em Riemannian join} of $S^{n-2}$ and $S^{1}$ with respect to
$\g$ and $\l$ and is denoted by $S^{n-2}*_{\gamma,\l} S^{1}$. 
\end{ede}

Notice that the metric $g$ extends to a continuous metric on $S^n$.
We will see below under which circumstances this extension is {\em smooth}.

\begin{ath} \label{t8} Let $N$ be a maximal leaf of the distribution
  $\xi ^\perp$ of $M_0$ and
  let $x\in N$ be the unique zero of the gradient conformal vector
  field $X=\n(|\xi|)$ on $N$. We then have 

$(i)$ There exists some positive number $l$, not depending on $N$, such
that the exponential map at $x$ maps
  diffeomorphically the open ball $B(0,l)$ in $T_xN$ onto $N$. 

$(ii)$ The submanifold $Z$ is connected, isometric to a round 
  sphere $S^{n-2}$. The closure of each integral leaf $N$ defined
  above is $\overline N=N\cup Z$. 

$(iii)$ $M$ is isometric to a Riemannian join $S^{n-2}*_{\gamma,\l}
S^{1}$, where $\g$ is $($up to a constant$)$ equal to the derivative of $\l$. 
\end{ath}

\bp  $(i)$ Consider the isometric action of $S^1$ on $M$ induced by
$\xi$. For $\theta\in
S^1$ denote by $N_\theta$ the image of $N$ through the action of
$\theta$ on $M$. Of course, $N_\theta$ is itself a maximal integral
leaf of $\xi ^\perp$. For every unit vector $V\in T_xN$, we
define 
$$l(x,V):=\sup\{t>0\ |\ \exp_x(rV)\in N,\ \forall\  r\le t\}.$$
Clearly, $l(x,V)$ is the distance along the geodesic $\exp_x(tV)$ from
$x$ to the first point on this geodesic where $X$ vanishes. Of course,
the exponential map on $N_\theta$ coincides (as long as it is 
defined) with the exponential map on $M$ since each $N_\theta$ is
totally geodesic. As noticed
before, the norm of $X$ along geodesics issued from $x$ only depends
on the parameter along the geodesic, therefore $l(x,V)$ is independent
of $V$ and can be denoted by $l(x)$. Since we have a transitive
isometric action on the $N_\theta$'s, $l(x)$ actually does not depend
on $x$ neither, and will be denoted by $l$. This proves that each $N_\theta$ is
equal to the image of the open ball $B(0,l)$ in
$T_{\theta(x)}N_\theta$ via the exponential $\exp_{\theta(x)}$. 

$(ii)$ Let us denote by $Z_\theta$ the set
$\overline{N_\theta}\setminus N_\theta$.
By the above, $Z_\theta$ is the image of the round sphere $S(0,l)$ in
$T_{\theta(x)}N_\theta\subset T_{\theta(x)}M$ via the exponential map
(on $M$) $\exp_{\theta(x)}$. In particular, each $Z_\theta$ is a
connected subset of $Z$,
diffeomorphic to $S^{n-2}$. Every element $\theta '\in S^1$ maps (by
continuity) $Z_\theta$ to $Z_{\theta '\theta}$ and on the other hand,
it preserves $Z$. We deduce that $Z_\theta=Z_{\theta '}$ for all
$\theta, \theta '\in S^1$ and since $Z$ is the union of all
$Z_\theta$, we obtain $Z=Z_\theta$. The other assertions are now
clear.

$(iii)$ This point is an implicit consequence of the local statements from
the previous sections. First, by Proposition \ref{prl} $M_0$ is
diffeomorphic to $N\times S^1$ with the warped product metric
$g=g_N+\lambda ^2 d\theta ^2$, where $\lambda$ is a function on $N$
whose gradient $X$ is a GCVF vanishing at $x$. From Proposition \ref{prds}
and $(i)$ above we see that $N\setminus \{x\}$ is diffeomorphic to $(0,l)\times
S^{n-2}$ with the metric $g_N=ds^2+ \gamma ^2(s){g_{S^{n-2}}}$. If we denote by
$S$ the orbit of $x$ under the $S^1$--action on $M$ defined by $\xi$,
this shows that $M_0\setminus S$ is diffeomorphic to $(0,l)\times S^{n-2}\times
S^1$ with the metric 
$$g=ds^2+ \gamma ^2(s){g_{S^{n-2}}} + \lambda ^2(s) d\theta ^2,$$
where $\l$ represents the norm of $\xi$ and
$X=\n\l=\l'(s)\frac{\pa}{\pa s}$. From (\ref{bg}) we get $|\l'|=|X|=c\gamma$.
Taking into account that $X$ does not vanish on $M_0\setminus S$, we see that
$\l'$ does not change sign on $(0,l)$, so $\g$ equals the derivative
of $\l$ up to some non--zero constant. Finally, the boundary
conditions (\ref{equiv1}) are easy to check: $\lim_{s\to 0} 
\gamma(s)=\frac1c |X_x|=0$, $\lim_{s\to l} \gamma(s)$ is equal to the radius of
the round $(n-2)$--sphere $Z$ and is thus positive, $\lim_{s\to 0} 
\l(s)=|\xi_x|>0$ and $\lim_{s\to l} 
\l(s)=0$ because $\xi$ vanishes on $Z$. 
\r

In order to obtain the classification we have to understand which of
the above Riemannian join metrics are actually smooth on the entire
manifold. For this we will use the following folkloric result

\begin{elem} \label{ext} Let $f:(0,\e)\to \RM^+$ be a smooth function
  such that $\lim_{t\to 0}f(t)=0$. The Riemannian metric
  $dt^2+f(t)g_{S^{n-1}}$ extends to a smooth metric at the
  singularity $t=0$ if and only if the function $\tilde
  f(t):=f(t^{\frac12})$ has a smooth extension at $t=0$ and $\tilde f'(0)=1$.
\end{elem}

Notice that the above condition on $f$ amounts to say that
$f(t)=t^2+t^4h(t^2)$ for some smooth germ $h$ around $0$.

\begin{ecor} Let $\g:(0,l)\to \RM^+$ be a smooth function satisfying 
$\lim_{s\to 0}\gamma(s)=0$ and 
$\lim_{s\to l}\gamma(s)>0$. For $c>0$ consider the function
\beq \label{l}\l(s):=c\int_s^l \g(t)dt.\eeq
The Riemannian join metric 
$$g=ds^2+\gamma ^2(s){g_{S^{n-2}}}+\l^2(s)d\theta ^2$$
defined on $(0,l)\times S^{n-2}\times S^1$
extends to a smooth metric on $S^n$ if and only if there exist two
smooth functions $a$ and $b$ defined on some interval $(-\e,\e)$ such that 
\beq \label{bou}
\g(t)=t(1+t^2a(t^2))\quad\hbox{and}\quad \g(l-t)=\frac{1}{c}+t^2b(t^2),\quad
\forall\  |t|<\e.
\eeq
\end{ecor}

\bp Since $\l(0)>0$, $g$ extends to a smooth metric at $s=0$ if and
only if $ds^2+\gamma ^2(s){g_{S^{n-2}}}$ extends smoothly at $s=0$. By
Lemma \ref{ext}, this is equivalent to $\g^2(t)=t^2+t^4h(t^2)$ for
some smooth $h$, so $\g(t)=t\sqrt{1+t^2h(t^2)}=t(1+t^2a(t^2))$ for
some smooth function $a$. Similarly, $g$ extends smoothly at $s=l$ if and only
if the same holds for $ds^2+\l^2(s)d\theta ^2$, which, by Lemma
\ref{ext} is equivalent to the existence of some smooth function
$d$ defined around $0$ such that $\l(l-t)=t+t^3d(t^2)$. Taking
(\ref{l}) into account, this is of course equivalent to the second
part of (\ref{bou}).
\r

Summarizing, we have

\begin{ath} \label{clas} Let $(M^n,g)$ be a compact Riemannian
  manifold carrying a 
  non--parallel Killing vector field $\xi$ whose covariant derivative
  is a twistor $2$--form. Then one of the following possibilities occurs:

1. $M$ is a space form of positive curvature and $\xi$ is any Killing
vector field on $M$.

2. $M$ is a Sasakian manifold and $\xi$ is the Sasakian vector field.

3. $M$ is a warped mapping torus $N_{\l,\f}$
$$M=(\RM\times N)/_{(t,x)\sim(t+1,\f(x))},\qquad g=\l ^2d\theta ^2+g_N,$$
 where $N$ is is a compact
$(n-1)$--dimensional Riemannian manifold carrying a GCVF with
primitive $\l$ $($cf. Remark \ref{re}$)$,
  $\f$ is an isometry of $N$ preserving $\l$ and $\xi =\frac{\pa}{\pa
    \theta}$.

4. $M$ is a Riemannian join $S^{n-2}*_{\gamma,\l} S^{1}$ 
with the metric $g=ds^2+\gamma ^2(s){g_{S^{n-2}}}+\l^2(s)d\theta ^2$
where $\g:(0,l)\to \RM^+$ is a smooth function satisfying the boundary
conditions $(\ref{bou})$, $\l$ is given by formula $(\ref{l})$, and
$\xi=\frac{\pa}{\pa \theta}$. 
\end{ath}

We end up these notes with some open problems related to the
classification above. One natural question is the following : which compact
Riemannian manifolds carry {\em twistor} $1$--forms $\xi$ with twistor exterior
derivative? To the author's knowledge, in all known examples $\xi$ is
either closed or co--closed. In the first case, the metric dual of
$\xi$ is a GCVF, so the manifold is described by Proposition
\ref{prds}. The second case just means that $\xi$ is Killing, and the
possible manifolds are described by Theorem \ref{clas}. 

More generally, one can address the question of classifying all
compact Riemannian manifolds $M^n$ carrying a Killing or twistor $p$--form
whose exterior derivative is a non--zero twistor form ($2\le
p\le n-2$). Besides the round spheres, the only known examples are
Sasakian manifolds (for odd $p$), nearly K\"ahler $6$--manifolds (for
$p=2$ and $p=3$) and nearly parallel $G_2$--manifolds (for $p=3$).


\vs

 \labelsep .5cm

\end{document}